\documentstyle[11pt,a4wide]{article}
\font \tenmsb=msbm10 scaled \magstep 1
\font \sevenmsb=msbm7 scaled \magstep 1
\font \fivemsb=msbm5 scaled \magstep 1
\newfam \msbfam
\textfont \msbfam=\tenmsb
\scriptfont \msbfam=\sevenmsb
\scriptscriptfont \msbfam=\fivemsb
\def \Bbb#1{\fam \msbfam \relax#1}

\font \teneufm=eufm10 scaled \magstep 1
\font \seveneufm=eufm7 scaled \magstep 1
\font \fiveeufm=eufm5 scaled \magstep 1
\newfam \eufmfam
\textfont \eufmfam=\teneufm
\scriptfont \eufmfam=\seveneufm
\scriptscriptfont \eufmfam=\fiveeufm
\def \frak#1{{\fam \eufmfam \relax#1}}

\def \unit{\hbox{$\hspace{.2em}{\Bbb I}\hspace{-.5em}\acute{}\hspace{.5em}$}}

\title{\bf ON FUNCTION THEORY IN QUANTUM DISC: A q-ANALOGUE OF BEREZIN
TRANSFORM}

\author{\sl D. Shklyarov \and \sl S. Sinel'shchikov
\and \sl L. Vaksman}

\date{\tt Institute for Low Temperature Physics \& Engineering\\
National Academy of Sciences of Ukraine}

\newtheorem{theorem}{Theorem}[section]
\newtheorem{lemma}[theorem]{Lemma}
\newtheorem{proposition}[theorem]{Proposition}
\newtheorem{corollary}[theorem]{Corollary}

\begin{document}

\maketitle

 Let $\alpha$ be a positive number.

 Section 1 of this work contains a study of Toeplitz-Bergman operators with
finite symbols in the quantum disc, and section 4 deals already with
Toeplitz-Bergman operators with bounded symbols. An alternate way of
producing Toeplitz-Bergman operators with polynomial symbols is described in
section 7 (lemma 7.2).

 Section 2 introduces a Berezin transform $B_{q,\alpha}$ for finite
functions in the quantum disc; the same is done in section 4 for bounded
functions. An alternate way of constructing a Berezin transform for a
polynomial function is described in sections 6, 7 (proposition 6.6 and lemma
7.2).

  An asymptotic expansion (3.2), (3.6) for a Berezin transform for a
finite function is obtained in section 3; a similar expansion (5.2)
for the case of a bounded function can be found in section 5. An application
of the latter result to formal series with polynomial coefficients in
section 8 affords the main result of this work (theorem 8.4).

 We use the background and notation used in \cite{SSV1, SSV2, SSV3}.

\bigskip

\section{Toeplitz-Bergman operators with finite symbols}

 Consider the covariant algebra ${\Bbb C}[z]_q$ (see \cite{SSV2}).
Algebraically it is isomorphic to the polynomial algebra ${\Bbb C}[z]$, and
the $U_q \frak{sl}_2$-action is determined by the relations
$$K^{\pm 1}z=q^{\pm 2}z,\qquad Fz=q^{1/2}.$$
We also follow \cite{SV} in using a covariant (left) ${\Bbb C}[z]_q$-module
with the generator $\unit$ and the relations
$$K^{\pm 1}\unit=q^{\pm(2 \alpha+1)}\unit,\qquad F \unit=0.$$
Denote this covariant module by ${\Bbb C}[z]_{q,\alpha}$

 Let $F_{q,\alpha}\subset{\rm End}({\Bbb C}[z]_{q,\alpha})$ be the covariant
algebra of linear operators $A:z^j \mapsto \sum \limits_{m \in {\Bbb
Z}_+}a_{mj}z^m$, $j \in{\Bbb Z}_+$, with finitely many nonzero matrix
elements $a_{mj}$. In virtue of this definition,
$F_{q,\alpha}\hookrightarrow{\Bbb C}[z]_{q,\alpha}\otimes{\Bbb
C}[z]_{q,\alpha}^*$

 Our immediate purpose is to construct a morphism of $U_q
\frak{sl}_2$-modules $D(U)_q \to F_{q,\alpha}$ which is normally called a
Toeplitz quantization.

 Remind the notation (see \cite{SSV1}):
$$\int \limits_{U_q}fd \nu_ \alpha=\frac{1-q^{4
\alpha}}{1-q^2}\int \limits_{U_q}f(1-zz^*)^{2 \alpha+1}d \nu,$$
$$(f_1,f_2)_{q,\alpha}=\int \limits_{U_q}f_2^*f_1d \nu_ \alpha.\eqno(1.1)$$
Form a completion of the linear space $D(U)_q$ with respect to the norm
$\|f\|_{q,\alpha}=(f,f)_{q,\alpha}^{1/2}$. It is easy to show that this
Hilbert space admits an embedding into $D(U)_q^\prime$ and is canonically
isomorphic to the space $L_{q,\alpha}^2$ defined in \cite{SSV1}.

 Let $P_{q,\alpha}$ be the orthogonal projection in $L_{q,\alpha}^2$ onto
the closure $H_{q,\alpha}^2$ of the subspace ${\Bbb C}[z]_{q,\alpha}\subset
L_{q,\alpha}^2$. Given $\stackrel{\circ}{f}\in D(U)_q$, we call the linear
operator
$$\widehat{f}:{\Bbb C}[z]_{q,\alpha}\to{\Bbb C}[z]_{q,\alpha};\qquad
\widehat{f}:\psi \mapsto P_{q,\alpha}(\stackrel{\circ}{f}\psi),\qquad \psi
\in{\Bbb C}[z]_{q,\alpha}$$
a Toeplitz-Bergman operator with the finite symbol $\stackrel{\circ}{f}$.
This is well defined, as one can see from

\medskip

\begin{proposition} With $\stackrel{\circ}{f}\in D(U)_q$, for all but
finitely many $m,j \in{\Bbb Z}_+$ the integral
$I_{m,j}=\displaystyle \int \limits_{U_q}z^{*m}\stackrel{\circ}{f}z^jd \nu_
\alpha$ is zero.
\end{proposition}

\smallskip

 {\bf Proof.} It was shown in \cite{SSV1} that for any
$\stackrel{\circ}{f}\in D(U)_q$ one has
$z^{*N}\stackrel{\circ}{f}=\stackrel{\circ}{f}z^N=0$ for some $N \in{\Bbb
N}$. Hence $I_{mj}=0$ if $\max(m,j)\ge N$. \hfill $\Box$

\medskip

 A straightforward consequence of proposition 1.1 is that the
Toeplitz-Bergman operator with a finite symbol belongs to the covariant
algebra $F_{q,\alpha}$.

\medskip

\begin{proposition} Toeplitz quantization $D(U)_q \to F_{q,\alpha}$,
$\stackrel{\circ}{f}\mapsto \widehat{f}$, is a morphism of $U_q
\frak{sl}_2$-modules.
\end{proposition}

\smallskip

 {\bf Proof.} One can deduce from the invariance of the scalar product in
${\Bbb C}[z]_{q,\alpha}$ and the covariance of the left ${\Bbb
C}[z]_q$-module ${\Bbb C}[z]_{q,\alpha}$ that the linear map
$$D(U)_q \otimes{\Bbb C}[z]_{q,\alpha}\to{\Bbb C}[z]_{q,\alpha};\qquad
\stackrel{\circ}{f}\otimes \psi \mapsto
P_{q,\alpha}(\stackrel{\circ}{f}\psi)$$
is a morphism of $U_q \frak{sl}_2$-modules. On the other hand, we need to
demonstrate that the linear map
$$D(U)_q \to{\Bbb C}[z]_{q,\alpha}\otimes{\Bbb C}[z]_{q,\alpha}^*,\qquad
\stackrel{\circ}{f}\mapsto \widehat{f}$$
(the tensor product here requires no completion due to proposition 1.1).
Observe that the two statements are equivalent to $U_q
\frak{sl}_2$-invariance of the same element of the corresponding completion
of the tensor product $F_{q,\alpha}\otimes D(U)_q^\prime$, which is
determined by the canonical isomorphisms ${\rm End}_{\Bbb C}(V_1,V_2)\simeq
V_2 \widehat{\otimes}V_1^*$, $(V_1 \otimes V_2)^*\simeq V_2^*
\widehat{\otimes}V_1^*$. \hfill $\Box$

\medskip

 Remind the notation ${\rm Fun}(U)_q={\rm Pol}({\Bbb C})_q+D(U)_q$. A very
important construction of \cite{SSV1} was the representation $T$ of ${\rm
Fun}(U)_q$ in the infinitely dimensional vector space $H$. A basis in $H$
was formed by the vectors $v_j=T(z^j)v_0$, $j \in{\Bbb Z}_+$ (see
\cite{SSV1}). $T$ provides a one-to-one map between the space of
finite functions $D(U)_q$ and the space of linear operators in $H$ whose
matrices in the basis $\{v_j \}_{j \in{\Bbb Z}_+}$ have finitely many
non-zero entries. For $j \in{\Bbb Z}_+$, let $f_j$ stand for such finite
function that $T(f_j)v_k=\delta_{jk}v_k$, $k \in{\Bbb Z}_+$.

 The relation $(1-zz^*)v_j=q^{2j}v_j$, $j \in{\Bbb Z}_+$, motivates the
following definition:
$$(1-zz^*)^\lambda \stackrel{def}{=}\sum_{n=0}^\infty q^{2n
\lambda}f_n,\qquad \lambda \in{\Bbb C}.$$
(The series converges in the topological space $D(U)_q^\prime$.)

 The work \cite{SSV2} presents an explicit form of the invariant integral in
the quantum disc. It is easy to show that for any finite function $f$
$$\int \limits_{U_q}fd \nu=(1-q^2)\,{\rm tr}\,T(f(1-zz^*)^{-1}).$$

\medskip

 {\bf Remark 1.3.} Let $\widehat{f}_0$ be the Toeplitz-Bergman operator with
symbol $f_0$. It follows from the relations $z^*f_0=f_0z=0$, $\displaystyle
\int \limits_{U_q}f_0d \nu=1-q^2$ that
$$\widehat{f}_0:z^j \mapsto \left \{\begin{array}{clc}1-q^{4 \alpha}&,&j=0\\
                                                     0&,&j \ne 0
                                    \end{array}\right..$$

\medskip \stepcounter{theorem}

 Now the relation (1.1), the trace properties and the definition of
$D(U)_q^\prime$ imply

\medskip

\begin{lemma}\hfill
\begin{enumerate}
\item For all $f \in D(U)_q$, $\lambda \in{\Bbb C}$
$$\int \limits_{U_q}f(z)(1-zz^*)^\lambda d \nu=\int
\limits_{U_q}(1-zz^*)^\lambda f(z)d \nu,$$
\item $$\int \limits_{U_q}f_1(z)f_2(z)(1-zz^*)d
\nu(z)=\int \limits_{U_q}f_2(z)f_1(z)(1-zz^*)d \nu(z)$$
for all $f_1(z)\in D(U)_q^\prime$, $f_2(z)\in D(U)_q$.
\end{enumerate}
\end{lemma}

\medskip

 The following proposition describes an integral representation for matrix
elements of Toeplitz-Bergman operator.

\medskip

\begin{proposition} Let $\stackrel{\circ}{f}\in D(U)_q$ and
$\widehat{f}:{\Bbb C}[z]_{q,\alpha}\to{\Bbb C}[z]_{q,\alpha}$;
$\widehat{f}:z^j \mapsto \sum \limits_{m \in{\Bbb Z}_+}\widehat{f}_{mj}z^m$
be a Toeplitz-Bergman operator with symbol $\stackrel{\circ}{f}$. Then
$$\widehat{f}_{mj}=\frac{1-q^{4
\alpha}}{1-q^2}\int \limits_{U_q}P_{z,mj}\stackrel{\circ}{f}(z)d
\nu(z),\eqno(1.2)$$
with
$$P_{z,mj}=\frac{(q^{4 \alpha+2};q^2)_m}{(q^2;q^2)_m}q^{2j}z^j(1-zz^*)^{2
\alpha+1}z^{*m}.\eqno(1.3)$$
\end{proposition}

\smallskip

 {\bf Proof.} Apply the relation \cite{KL, SSV1}
$$(z^m,z^l)_{q,\alpha}=\frac{(q^2;q^2)_m}{(q^{4
\alpha+2};q^2)_m}\delta_{ml},\qquad m,l \in{\Bbb Z}_+$$
to get
$$\widehat{f}_{mj}=
\frac{(\stackrel{\circ}{f}z^j,z^m)_{q,\alpha}}{(z^m,z^m)_{q,\alpha}}=
\frac{(q^{4 \alpha+2};q^2)_m}{(q^2;q^2)_m}\int
\limits_{U_q}z^{*m}\stackrel{\circ}{f}z^jd \nu_ \alpha=$$
$$=\frac{1-q^{4 \alpha}}{1-q^2}\cdot \frac{(q^{4
\alpha+2};q^2)_m}{(q^2;q^2)_m}\int
\limits_{U_q}z^{*m}\stackrel{\circ}{f}z^j(1-zz^*)^{2 \alpha+1}d \nu.$$
Hence, by lemma 1.4,
$$\widehat{f}_{mj}=\frac{1-q^{4 \alpha}}{1-q^2}\cdot \frac{(q^{4
\alpha+2};q^2)_m}{(q^2;q^2)_m}\int \limits_{U_q}(1-zz^*)z^j(1-zz^*)^{2
\alpha}z^{*m}\stackrel{\circ}{f}d \nu.\eqno(1.4)$$
It remains to apply the relation
$$(1-zz^*)z=q^2z(1-zz^*).\eqno \Box$$

\medskip

 {\bf Remark 1.6.} The matrix $P_z=(P_{z,mj})_{m,j \in{\Bbb Z}_+}$ is a
q-analogue for the matrix of a one-dimensional orthogonal projection onto
the subspace generated by the vector $k_z$ from an overfull system (see
\cite{B}). A q-analogue of the overfull system itself is presented in the
Appendix.

\medskip

 {\bf Remark 1.7.} It follows from proposition 1.2 and the relation $U_q
\frak{sl}_2 \cdot \widehat{f}_0=F_{q,\alpha}$ to be proved later on (see
proposition 6.4) that the map $D(U)_q \to F_{q,\alpha}$,
$\stackrel{\circ}{f}\mapsto \widehat{f}$, given by Toeplitz quantization is
onto.

\bigskip

\section{Berezin transform: finite functions}

 Consider a $U_q \frak{sl}_2$-module $V$ and the covariant algebra ${\rm
End}_{\Bbb C}(V)_f \simeq V \otimes V^*$. There is a well known (see
\cite{CP, SSV2}) formula for an invariant integral
$${\rm tr}_q:{\rm End}_{\Bbb C}(V)_f \to{\Bbb C},\qquad {\rm tr}_q:A
\mapsto{\rm tr}(A \cdot K^{-1}).$$
In the case $V={\Bbb C}[z]_{q,\alpha}$ and $A:z^j \mapsto \sum \limits_{m
\in{\Bbb Z}_+}a_{mj}z^m$ being an element of the covariant algebra
$F_{q,\alpha}\subset{\rm End}_{\Bbb C}({\Bbb C}[z]_{q,\alpha})_f$, one has
${\rm tr}_q(A)=\sum \limits_{k \in{\Bbb Z}_+}a_{kk}q^{-2k}$.

 Given a linear operator $\widehat{f}\in F_{q,\alpha}$, a distribution $f
\in D(U)_q^\prime$ is said to be a symbol of $\widehat{f}$ if for all
$\stackrel{\circ}{\psi}\in D(U)_q$
$$\int \limits_{U_q}f \cdot \stackrel{\circ}{\psi}d \nu=\frac{1-q^2}{1-q^{4
\alpha}}{\rm tr}_q(\widehat{f}\widehat{\psi}).\eqno(2.1)$$
(Here $\widehat{\psi}$ is the Toeplitz-Bergman operator with symbol
$\stackrel{\circ}{\psi}$.)

 This definition is a q-analogue of the Berezin's definition, as one can
observe from relation (3.15) from \cite{B}.

\medskip

\begin{proposition} The covariant symbol of a linear operator
$\widehat{f}:z^j \mapsto \sum \limits_{m \in{\Bbb Z}_+}\widehat{f}_{mj}z^m$,
$j \in{\Bbb Z}_+$, from the algebra $F_{q,\alpha}$, is given by
$$f={\rm tr}_q(\widehat{f}\cdot P_z)=\sum_{j,m \in{\Bbb
Z}_+}\widehat{f}_{jm}P_{z,mj}q^{-2j}.\eqno(2.2)$$
\end{proposition}

\smallskip

 {\bf Proof.} By a virtue of (1.2)
$${\rm tr}_q(\widehat{f}\widehat{\psi})=\sum_{j,m \in{\Bbb
Z}_+}\widehat{f}_{jm}\widehat{\psi}_{mj}q^{-2j}=\sum_{j,m \in{\Bbb
Z}_+}\widehat{f}_{jm}\frac{1-q^{4 \alpha}}{1-q^2}\int
\limits_{U_q}P_{z,mj}\stackrel{\circ}{\psi}d \nu(z)\cdot q^{-2j}=$$
$$=\int \limits_{U_q}\left(\frac{1-q^{4 \alpha}}{1-q^2}\sum_{j,m \in{\Bbb
Z}_+}\widehat{f}_{jm}P_{z,mj}\cdot q^{-2j}\right)\stackrel{\circ}{\psi}d
\nu(z).\eqno \Box$$

\medskip

 Note that the integral representation (2.2) is a q-analogue of the relation
(3.4) in \cite{B}.

 On can deduce from the covariance of algebras $D(U)_q$, $F_{q,\alpha}$, the
invariance of the integrals $\nu:D(U)_q \to{\Bbb C}$, ${\rm
tr}_q:F_{q,\alpha}\to{\Bbb C}$, the "integration in parts" formula
\cite[proposition 2.1]{SSV2}, and proposition 1.2 the following

\medskip

\begin{proposition} The linear map $F_{q,\alpha}\to D(U)_q^\prime$,
$\widehat{f}\mapsto f$, which takes a linear operator to its covariant
symbol, is a morphism of $U_q \frak{sl}_2$-modules.
\end{proposition}

\medskip

 As in \cite{UU}, we call the covariant symbol $f$ for the Toeplitz-Bergman
operator $\widehat{f}$ with symbol $\stackrel{\circ}{f}\in D(U)_q$ a Berezin
transform of the function $\stackrel{\circ}{f}$. The associated
transform map will be denoted by $B_{q,\alpha}$:
$$B_{q,\alpha}:D(U)_q \to D(U)_q^\prime;\qquad
B_{q,\alpha}:\stackrel{\circ}{f}\:\mapsto f.$$

 Propositions 1.2 and 2.2 imply

\medskip

\begin{proposition} The Berezin transform is a morphism of $U_q
\frak{sl}_2$-modules.
\end{proposition}

\medskip

 {\bf Example 2.4.} Let $\widehat{f}\in F_{q,\alpha}$ be given by
$\widehat{f}z^j=\left \{{\textstyle 1,j=0 \atop \textstyle 0,j \ne
0}\right.$. Then one has $f=(1-zz^*)^{2 \alpha+1}$. Hence,
$B_{q,\alpha}f_0=(1-q^{4 \alpha})(1-zz^*)^{2 \alpha+1}$ since
$\widehat{f_0}:z^j \mapsto \left \{\begin{array}{clc}1-q^{4 \alpha}&,&j=0\\
                                                     0&,&j \ne 0
                                    \end{array}\right.$ (see Example 1.3).

\medskip \stepcounter{theorem}

 To conclude, we prove that Berezin transform is an integral operator,
and find its kernel. In this way, a q-analogue of the relation (4.8) from
\cite{B} is to be obtained.

\medskip

\begin{proposition} For all $\stackrel{\circ}{f}\in D(U)_q$,
$$(B_{q,\alpha}\stackrel{\circ}{f})(z)=\int
\limits_{U_q}b_{q,\alpha}(z,\zeta)\stackrel{\circ}{f}(\zeta)d \nu(\zeta),$$
with $b_{q,\alpha}\in D(U \times U)_q^\prime$ being given by
$$b_{q,\alpha}(z,\zeta)=\frac{1-q^{4 \alpha}}{1-q^2}(1-zz^*)^{2
\alpha+1}(1-\zeta \zeta^*)^{2 \alpha+1}\{(q^2z^*\zeta;q^2)_{-(2 \alpha+1)}
\cdot(z \zeta^*;q^2)_{-(2 \alpha+1)}\}.$$
(See \cite{SSV3} for the definition of $\{.,.\}$.)
\end{proposition}

\smallskip

 {\bf Proof.} Consider the linear operator
$$\widetilde{B}_{q,\alpha}:D(U)_q \to D(U)_q^\prime;\qquad
\widetilde{B}_{q,\alpha}:\stackrel{\circ}{f}\mapsto \int
\limits_{U_q}b_{q,\alpha}(z,\zeta)\stackrel{\circ}{f}(\zeta)d \nu(\zeta).$$
Its kernel coincides up to a constant multiple to the invariant kernel
$k_{22}^{-(2 \alpha+1)}\cdot k_{11}^{-(2 \alpha+1)}$ (see \cite{SSV3}).
Hence, $\widetilde{B}_{q,\alpha}$ is a morphism of $U_q
\frak{sl}_2$-modules by \cite[proposition 4.5]{SSV2}. Note that
$B_{q,\alpha}$ possesses the same property. It was shown in \cite{SSV2} that
$f_0 \in D(U)_q$ generates the $U_q \frak{sl}_2$-module $D(U)_q$. In this
context, the desired equality $\widetilde{B}_{q,\alpha}=B_{q,\alpha}$
becomes a consequence of $\widetilde{B}_{q,\alpha}f_0=(1-q^{4
\alpha})(1-zz^*)^{2 \alpha+1}=B_{q,\alpha}f_0$. \hfill $\Box$

\bigskip

\section{Berezin transform and Laplace-Beltrami operator}

 The following lemma is deduced from the relation
$$(1-zz^*)^\lambda{=}\sum_{n=0}^\infty q^{2n \lambda}f_n,\qquad \lambda
\in{\Bbb C}$$
and (1.3):

\medskip

\begin{lemma} For all $m,j \in{\Bbb Z}_+$ the following decomposition is
valid in $D(U)_q^\prime$:
$$P_{z,mj}=\frac{(q^{4 \alpha+2};q^2)_m}{(q^2;q^2)_m}\sum_{n=0}^\infty q^{4
\alpha n}\cdot P_{z,mj}^{(n)},\eqno(3.1)$$
with $P_{z,mj}^{(n)}=q^{2(j+n)}z^j \cdot f_n \cdot z^{*m} \in D(U)_q$.
\end{lemma}

\medskip

 Let $\stackrel{\circ}{f},\psi \in D(U)_q$. Consider the integral
$\displaystyle \int \limits_{U_q}\psi^*\cdot
B_{q,\alpha}\stackrel{\circ}{f}d \nu$ as a function of $t=q^{4 \alpha}$. Now
proposition 2.1 and lemma 3.1 imply the analyticity of this function as $t
\in[0,1)$. Hence, one has

\medskip

\begin{proposition} There exists a unique sequence of $U_q
\frak{sl}_2$-module morphisms $B_q^{(n)}:D(U)_q \to D(U)_q^\prime$, $n
\in{\Bbb Z}_+$, such that for all $\stackrel{\circ}{f}\in D(U)_q$
$$B_{q,\alpha}\stackrel{\circ}{f}=\sum_{n=0}^\infty q^{4 \alpha
n}B_q^{(n)}\stackrel{\circ}{f}.\eqno(3.2)$$
\end{proposition}

\medskip

 Our purpose is to prove that the linear operators $B_q^{(n)}$ are
polynomials of Laplace-Beltrami operator in the quantum disc.

 Let
$$p_j(t)= \sum_{k=0}^j \frac{(q^{-2j};q^2)_k}{(q^2;q^2)^2_k}q^{2k}
\cdot
\prod_{i=0}^{k-1}\left(1-q^{2i}\left((1-q^2)^2t+1+q^2
\right)+q^{4i+2}\right).\eqno(3.3)$$

\medskip

\begin{lemma} $ p_j(\Box)f_0=q^{2j}\cdot f_j$ for all $j \in{\Bbb Z}_+$.
\end{lemma}

\smallskip

 {\bf Proof.} Remind \cite{SSV1} that for all $l \in{\Bbb C}$ the basic
hypergeometric series
$$\varphi_l=_3\!\!\Phi_2 \left[{(1-zz^*)^{-1},q^{-2l},q^{2(l+1)};q^2;q^2
\atop q^2,0}\right]$$
converge in $D(U)_q^\prime$, and
$$\Box \varphi_l=-\frac{(1-q^{-2l})(1-q^{2l+2})}{(1-q^2)^2}\varphi_l.$$

 By a virtue of \cite[\S 6]{SSV1}, it suffices to show that for all $l
\in{\Bbb C}$
$$q^{-2j}\cdot\int \limits_{U_q}\varphi_l^*\cdot p_j(\Box)f_0d \nu=\int
\limits_{U_q}\varphi_l^*f_jd \nu.\eqno(3.4)$$
After substituting $l$ by $\overline{l}$ we find out that (3.4) is
equivalent to
$$p_j \left(-\frac{(1-q^{-2l})(1-q^{2(l+1)})}{(1-q^2)^2}\right)=\,_3
\Phi_2 \left[{q^{-2j},q^{-2l},q^{2(l+1)};q^2;q^2 \atop q^2,0}\right].$$

 Prove this relation. By the definition of $_3 \Phi_2$ one has
$$_3 \Phi_2 \left[{q^{-2j},q^{-2l},q^{2(l+1)};q^2;q^2 \atop q^2,0}\right]=$$
$$=\sum_{k=0}^j \frac{(q^{-2j};q^2)_k}{(q^2;q^2)^2_k}\cdot
\prod_{i=0}^{k-1}\left((1-q^{2i}\cdot q^{2(l+1)})(1-q^{2i}\cdot
q^{-2l})\right)\cdot q^{2k}=$$
$$=\sum_{k=0}^j \frac{(q^{-2j};q^2)_k}{(q^2;q^2)^2_k}\cdot
\prod_{i=0}^{k-1}\left((1+q^{2i}\cdot u+q^{4i+2}\right)\cdot q^{2k}$$
with $u=-q^{2l+2}-q^{-2l}$. It remains to prove that
$$p_j \left(-\frac{1+q^2+u}{(1-q^2)^2}\right)=\sum_{k=0}^j
\frac{(q^{-2j};q^2)_k}{(q^2;q^2)^2_k}\cdot
\prod_{i=0}^{k-1}\left((1+q^{2i}\cdot u+q^{4i+2}\right)\cdot q^{2k}.$$
For that, it suffices to exclude $u$ by a substitution $u=-
(1-q^2)^2t-1-q^2.$ \hfill $\Box$

\medskip

 The next statement refines essentially proposition 3.2.

\medskip

\begin{proposition} For all $\stackrel{\circ}{f}\in D(U)_q$ the following
expansion in $D(U)_q^\prime$ is valid:
$$B_{q,\alpha}\stackrel{\circ}{f}=(1-q^{4 \alpha})\sum_{j \in{\Bbb
Z}_+}q^{4 \alpha \cdot j}\cdot p_j(\Box)\stackrel{\circ}{f}.\eqno(3.5)$$
\end{proposition}

\smallskip

 {\bf Proof.} One has the relation $B_{q,\alpha}f_0=(1-q^{4
\alpha})\displaystyle \sum \limits_{k \in{\Bbb Z}_+}q^{(4 \alpha+2)k}f_k$
(see example 2.4). Hence, in the special case $\stackrel{\circ}{f}\,=f_0$ our
statement follows from lemma 3.3.  It remains to take into account that
$f_0$ generates the $U_q \frak{sl}_2$-module $D(U)_q$, and the operators
$B_{q,\alpha}$, $\Box$ are morphisms of $U_q \frak{sl}_2$-modules (see
\cite[proposition 2.1]{SSV1}).  \hfill $\Box$

\medskip

\begin{corollary}
$$B_q^{(n)}=\left \{\begin{array}{lcl}I &,& n=0 \\
                                      p_n(\Box)-p_{n-1}(\Box) &,&
                                      n \in{\Bbb
N}\end{array}\right.\eqno(3.6)$$
\end{corollary}

\bigskip

\section{Toeplitz-Bergman operators with bounded symbols}

 It is very well known \cite{KL, NN} that the $^*$-algebra ${\rm Pol}({\Bbb
C})_q$ has a unique up to unitary equivalence faithful irreducible
representation. As it was described in \cite{SSV1}, this representation $T$
lives in a Hilbert space $\overline{H}$ constructed as a completion of the
pre-Hilbert space $H$. Let $L(\overline{H})$ be the algebra of all bounded
operators in $H$ and $H^\prime$ the vector space of all bounded antilinear
functionals on $H$. One has
$${\rm End}_{\Bbb C}(H)\subset L(\overline{H})\subset{\rm Hom}_{\Bbb
C}(H,H^\prime).$$
It was demonstrated in \cite{SSV1} that the map $T:{\rm Pol}({\Bbb C})_q
\hookrightarrow L(\overline{H})$ is extendable by a continuity up to the
isomorphism $T:D(U)_q^\prime {{\atop \textstyle \to} \atop {\textstyle
\approx \atop}}{\rm Hom}_{\Bbb C}(H,H^\prime)$.

 We call a distribution $f \in D(U)_q^\prime$ bounded if $T(f)\in
L(\overline{H})$. Impose the notation
$$L_q^\infty=\{f \in D(U)_q^\prime|\:T(f) \in L(\overline{H})\},\qquad \|f
\|_\infty=\|T(f)\|.$$
(It is easy to show that the algebra $L^\infty$ defined in this way is
isomorphic to the enveloping von Neumann algebra of the $C^*$-algebra of
continuous functions in the quantum disk, which was considered, in
particular, in \cite{NN}).

 Consider the subspaces
$${\Bbb C}[z]_{q,\infty}=\{f \in D(U)_q|\:f \cdot z=0 \},$$
$$H_{q,\infty}^2=\{f \in L^2(U)_q|\:f \cdot z=0 \},$$
$${\Bbb C}[[z]]_{q,\infty}=\{f \in D(U)_q^\prime|\:f \cdot z=0 \}.$$
It follows from \cite[proposition 3.3]{SSV2} that
$${\Bbb C}[z]_{q,\infty}={\Bbb C}[z]\cdot f_0,\qquad{\Bbb
C}[[z]]_{q,\infty}={\Bbb C}[[z]]\cdot f_0,$$
and hence
$$H \simeq{\Bbb C}[z]_{q,\infty},\qquad \overline{H}\simeq
H_{q,\infty}^2,\qquad H^\prime \simeq{\Bbb C}[[z]]_{q,\infty}.$$

 $T$ is unitarily equivalent to the representation $\widehat{T}$ of ${\rm
Pol}({\Bbb C})_q$ in $H_{q,\infty}^2$ given by
$$\widehat{T}:\psi \mapsto f \cdot \psi;\qquad f \in {\rm Pol}({\Bbb
C})_q,\;\psi \in H_{q,\infty}^2 \subset D(U)_q^\prime.$$
Thus, a distribution $f \in D(U)_q^\prime$ is bounded iff the linear
operator $\widehat{T}(f)$ is in $L(H_{q,\infty}^2)$; in this case $\|f
\|_\infty=\|\widehat{T}(f)\|_\infty$.

 The following proposition justifies the use of the symbol $\infty$ in the
notation for the vector spaces ${\Bbb C}[z]_{q,\infty}$, $H_{q,\infty}^2$,
${\Bbb C}[[z]]_{q,\infty}$.

\medskip

\begin{proposition} For any polynomial $\psi \in{\Bbb C}[z]_q$
$$\lim_{\alpha \to \infty}(\psi,\psi)_{q,\alpha}={1 \over 1-q^2}(\psi
f_0,\psi f_0).$$
\end{proposition}

\smallskip

 {\bf Proof.}
$$(\psi,\psi)_{q,\infty}\stackrel{\rm def}{=}\lim_{\alpha \to \infty}{1-q^2
\over 1-q^{4 \alpha}}(\psi,\psi)_{q,\alpha}=\lim_{\alpha \to \infty}\int
\limits_{U_q}\psi^*\psi \sum_{n=0}^\infty q^{4n \alpha}f_nd \nu=$$
$$=\int \limits_{U_q}\psi^*\psi f_0d \nu=(\psi f_0,\psi f_0).\eqno \Box$$

\medskip

 The following remark will not be used in the sequel. Proposition 4.1 allows
one to prove that the covariant algebra $D(U)_q$ is isomorphic to a "limit
$F_{q,\infty}$ of covariant algebras $F_{q,\alpha}$ as $\alpha \to \infty$".
This leads to an alternate scheme of producing the covariant algebra
$D(U)_q$ of finite functions in the quantum disk. Under this scheme, at
the first step a unitarizable Harish-Chandra module $V_ \alpha$ with lowest
weight $\alpha>0$ and the covariant algebras $V_ \alpha \otimes V_ \alpha^*
\hookrightarrow{\rm End}_{\Bbb C}(V_ \alpha)$ are constructed. The second
step is in "passage to the limit" $\lim \limits_{\alpha \to+\infty}V_ \alpha
\otimes V_ \alpha^*$ which is to be declared the algebra of finite
functions in the quantum disk.

 Finally, impose the notation
$$\overline{F}_{q,\infty}\stackrel{\rm def}{=}{\rm End}_{\Bbb C}({\Bbb
C}[z]_{q,\infty},{\Bbb C}[[z]]_{q,\infty}).$$
It follows from the definitions that the representation $\widehat{T}$ is
extendable up to a bijection $\widehat{T}:D(U)^\prime{{\atop \textstyle \to}
\atop{\textstyle \sim \atop}}\overline{F}_{q,\infty}$.

 It should be noted that ${\rm Pol}({\Bbb C})_q \subset L_q^\infty$. This
can be deduced, for example, from the fact that the representation
$\widehat{T}$ of ${\rm Pol}({\Bbb C})_q$ in the pre-Hilbert space ${\Bbb
C}[z]_{q,\infty}$ is a $^*$-representation of this algebra. Hence,
$I-\widehat{T}(z)\widehat{T}(z^*)\ge 0$,
$\|\widehat{T}(z)\|=\|\widehat{T}(z^*)\|=1$.

 Let $A$ be a compact linear operator in a Hilbert space and
$|A|\stackrel{\rm def}{=}(A^*A)^{1/2}$. Consider the sequence of eigenvalues
of $|A|$, with their multiplicities being taken into account:
$$s_1(A)\ge s_2(A)\ge \ldots.$$
The numbers $s_p(A)$, $p \in{\Bbb N}$, are called s-values of $A$.

 Remind the notation $S_ \infty$ for the ideal of all compact operators in
a Hilbert space, together with the notation
$$\|A \|_p=\left(\sum_{n \in{\Bbb N}}s_n(A)^p\right)^{1/p},\qquad S_p=\left
\{A \in S_ \infty|\:\|A \|_p<\infty \right \},\qquad p>0,$$
for the normed ideals of von Neumann-Schatten (see \cite{DS}).

\medskip

\begin{lemma} For any function $\psi \in D(U)_q$
$$\|\psi \|=(1-q^2)^{1/2}\cdot \|\widehat{T}(\psi(1-zz^*)^{-1/2})\|_2,$$
with $\|\psi \|=\left(\displaystyle \int \limits_{U_q}\psi^*\psi d \nu
\right)^{1/2}$.
\end{lemma}

\smallskip

 {\bf Proof.} It follows from (1.1) and the well known tracial properties of
an operator $A \in S_1$ that
$$\|\psi \|^2=(1-q^2){\rm
tr}\,\widehat{T}(\psi^*\psi(1-zz^*)^{-1})=(1-q^2){\rm
tr}\,\widehat{T}((1-zz^*)^{-1/2}\psi^*\psi(1-zz^*)^{-1/2}).\eqno \Box$$

\medskip

\begin{corollary} Let $\stackrel{\circ}{f}\in L_q^\infty$, $\psi \in
D(U)_q$, then $\stackrel{\circ}{f}\psi \in L^2(d \nu)_q$ and
$\|\stackrel{\circ}{f}\psi \|\le \|\stackrel{\circ}{f}\|_ \infty \cdot \|\psi
\|$.
\end{corollary}

\smallskip

 {\bf Proof.}
$$\|\stackrel{\circ}{f}\psi
\|=(1-q^2)^{1/2}\|\widehat{T}(\stackrel{\circ}{f})\widehat{T}
(\psi(1-zz^*)^{-1/2})\|_2 \le$$
$$\le(1-q^2)^{1/2}\|\widehat{T}(\stackrel{\circ}{f})\|\cdot
\|\widehat{T}(\psi(1-zz^*)^{-1/2})\|_2=\|\stackrel{\circ}{f}\|_ \infty \cdot
\|\psi \|.\eqno \Box$$

\medskip

 It follows from the boundedness of the multiplication operator by a bounded
function $\stackrel{\circ}{f}$:
$$D(U)_q \to L^2(d \nu)_q,\qquad \psi \mapsto \stackrel{\circ}{f}\psi$$
that it admits an extension by a continuity onto the entire space $L^2(d
\nu)_q$. This allows one to define a Toeplitz-Bergman operator $\widehat{f}$
with symbol $\stackrel{\circ}{f}\in L_q^\infty$:
$$\widehat{f}:H_{q,\alpha}^2 \to H_{q,\alpha}^2;\qquad \widehat{f}:\psi
\mapsto P_{q,\alpha}(\stackrel{\circ}{f}\psi).$$
By a virtue of corollary 4.3 one has
$$\|\widehat{f}\|\le \|\stackrel{\circ}{f}\|_ \infty,\eqno(4.1)$$
with $\|\widehat{f}\|$ being the norm of the operator $\widehat{f}$ in
$H_{q,\alpha}^2$. Thus we get a norm decreasing linear map $L_q^\infty \to
L(H_{q,\alpha}^2)$, $\stackrel{\circ}{f}\mapsto \widehat{f}$. This
definition generalizes that of a Toeplitz-Bergman operator with finite
symbol (see section 2).

\bigskip

\section{Berezin transform: bounded functions}

 The $U_q \frak{sl}_2$-module ${\Bbb C}[z]_{q,\alpha}$ is formed by
polynomials $\psi=\sum \limits_{i \in{\Bbb Z}_+}a_i(\psi)z^i$. Consider a
completion ${\Bbb C}[[z]]_{q,\alpha}$ of this vector space in the topology
of coefficientwise convergence, and impose the notation
$\overline{F}_{q,\alpha}\stackrel{\rm def}{=}{\rm Hom}_{\Bbb C}({\Bbb
C}[z]_{q,\alpha},{\Bbb C}[[z]]_{q,\alpha})$ for the corresponding
completion of $F_{q,\alpha}$. Equip $\overline{F}_{q,\alpha}$ with the
topology of pointwise (strong) convergence:
$$\lim_{n \to \infty}A_n=A \quad \Leftrightarrow \quad \forall \psi
\in{\Bbb C}[z]_{q,\alpha}\;\lim_{n \to \infty}A_n \psi=A \psi.$$
Evidently, ${\Bbb C}[z]_{q,\alpha}\subset H_{q,\alpha}^2 \subset{\Bbb
C}[[z]]_{q,\alpha}$, and so
$$F_{q,\alpha}\subset L(H_{q,\alpha}^2)\subset \overline{F}_{q,\alpha}.$$
The representation operators of $E$, $F$, $K^{\pm 1}$ in ${\Bbb
C}[z]_{q,\alpha}$ have degrees $+1$, $-1$, 0 respectively. Hence they are
extendable by a continuity from ${\Bbb C}[z]_{q,\alpha}$ onto ${\Bbb
C}[[z]]_{q,\alpha}$, and from $F_{q,\alpha}$ onto $\overline{F}_{q,\alpha}$.

 Of course, $\overline{F}_{q,\alpha}$ is a covariant bimodule over the
covariant algebra $F_{q,\alpha}$. It is easy to show that the linear
functional
$$F_{q,\alpha}\otimes F_{q,\alpha}\to{\Bbb C},\qquad \widehat{f}\otimes
\widehat{\psi}\mapsto {\rm tr}_q(\widehat{f}\widehat{\psi})$$
is extendable by a continuity up to a morphism of $U_q \frak{sl}_2$-modules
$\overline{F}_{q,\alpha}\otimes F_{q,\alpha}\to{\Bbb C}$.

 Define a covariant symbol $f \in D(U)_q^\prime$ of a linear operator
$\widehat{f}\in \overline{F}_{q,\alpha}$ by (2.1). The map
$\overline{F}_{q,\alpha}\to D(U)_q^\prime$ arising this way is a $U_q
\frak{sl}_2$-module morphism.

 In the following proposition we use notation $\widehat{f}$ for a linear
operator without assuming it to be a Toeplitz-Bergman operator.

\medskip

\begin{proposition} Let $\widehat{f}$ be a linear operator
$$\widehat{f}:{\Bbb C}[z]_{q,\alpha}\to{\Bbb C}[[z]]_{q,\alpha},\qquad
\widehat{f}:z^j \mapsto \sum_{m \in{\Bbb Z}_+}\widehat{f}_{mj}z^m,\quad j
\in{\Bbb Z}_+.$$
The series $\displaystyle \sum \limits_{j,m \in{\Bbb
Z}_+}\widehat{f}_{jm}P_{z,mj}q^{-2j}$ converges in $D(U)_q^\prime$ to the
covariant symbol of $\widehat{f}$.  \end{proposition}

\smallskip

 {\bf Proof.} It follows from the results of section 1 that for any
$\stackrel{\circ}{\psi}\in D(U)_q$ all but finitely many of integrals
$\displaystyle \int \limits_{U_q}P_{z,mj}\stackrel{\circ}{\psi}(z)d \nu(z)$
are zero. This allows one to reproduce literally the argument used in the
proof of proposition 2.1. \hfill $\Box$

\medskip

 Let $\stackrel{\circ}{f}\in L_q^\infty$, and $\widehat{f}\in
L(H_{q,\alpha}^2)\subset \overline{F}_{q,\alpha}$ be the Toeplitz-Bergman
operator with symbol $\stackrel{\circ}{f}$. We follow \cite{UU} in using the
term "Berezin transform of the function $\stackrel{\circ}{f}$" for the
covariant symbol of the linear operator $\widehat{f}$.

 Our purpose is to decompose the operator-function $B_{q,\alpha}:L_q^\infty
\to D(U)_q^\prime$ into series in powers of $t=q^{4 \alpha}$ (cf. (3.5)).

 One can use again the argument of proposition 1.5 to get (1.4) for all
bounded symbols $\stackrel{\circ}{f}\in L_q^\infty$. An application of (1.1)
and the fact that $\widehat{T}((1-zz^*)^{2 \alpha})$ is a trace class
operator for all $\alpha>0$, yields also

\medskip

\begin{proposition} For all $\stackrel{\circ}{f}\in L_q^\infty$ $m,j
\in{\Bbb Z}_+$
$$\widehat{f}_{mj}=
\frac{(q^{4 \alpha+2};q^2)_m}{(q^2;q^2)_m}\cdot(1-q^{4
\alpha})\cdot{\rm tr}\left(\widehat{T}(z^j(1-zz^*)^{2
\alpha}z^{*m})\widehat{T}(\stackrel{\circ}{f})\right).$$
\end{proposition}

\medskip

 Let $\Theta$ be the vector space of holomorphic functions in the unit disc
with values in the Banach algebra $S_1$ of trace class operators in
$\overline{H}$. (Each function $Q(t)$ from $\Theta$ admits a expansion
into the power series $Q(t)=\sum \limits_{n \in{\Bbb Z}_+}t^n \cdot Q^{(n)}$
with $\lim \limits_{n \to \infty}\|Q^{(n)}\|_1^{1/n}\le 1$.)

\medskip

\begin{proposition} For all $j,m \in{\Bbb Z}_+$
$$\sum_{n \in{\Bbb Z}_+}t^n \cdot \widehat{T}(z^j \cdot f_n \cdot z^{*m})
\in \Theta.$$
\end{proposition}

\smallskip

 {\bf Proof.} Remind that
$\widehat{T}(z)\widehat{T}(z^*)=1-\displaystyle \sum \limits_{n \in{\Bbb
Z}_+}q^{2n}\cdot \widehat{T}(f_n)$, and that $\widehat{T}(f_n)$ are
one-dimensional projections, $n \in{\Bbb Z}_+$. Hence
$\|\widehat{T}(z)\|=\|\widehat{T}(z^*)\|=\|\widehat{T}(f_n)\|_1=1$. Finally,
$$\|\widehat{T}(z^jf_nz^{*m})\|_1 \le\|\widehat{T}(z)\|^j \cdot
\|\widehat{T}(f_n)\|_1 \cdot \|\widehat{T}(z^*)\|^m=1.\eqno \Box$$

\medskip

 Propositions 5.2, 5.3 and the definition of Berezin transform imply

\medskip

\begin{corollary} Let $\psi \in D(U)_q$. There exists a unique function $Q_
\psi(t)\in \Theta$ such that
$$\int \limits_{U_q}(B_{q,\alpha}\stackrel{\circ}{f})\psi d \nu={\rm
tr}\left(\widehat{T}(\stackrel{\circ}{f})Q_ \psi(q^{4
\alpha})\right)\eqno(5.1)$$
for all $\stackrel{\circ}{f}\in L_q^\infty$.
\end{corollary}

\smallskip

 {\bf Proof.} The uniqueness of $Q_ \psi(t)$ is evident. In fact, given such
$A \in L(\overline{H})$ that for all $\stackrel{\circ}{f}\in D(U)_q$ one has
${\rm tr}\left(\widehat{T}(\stackrel{\circ}{f})A \right)=0$, then surely
$A=0$. The existence of $Q_ \psi \in \Theta$ follows from propositions 5.2,
5.3 and the definition of Berezin transform. \hfill $\Box$

\medskip

 The coefficients of the Taylor series for the holomorphic function $Q_
\psi(t)$ at $t=0$ are trace class operators. One can use (3.5) and (1.1) to
express those coefficients via the operators $\widehat{T}(p_j(\Box)\psi)$,
$j \in{\Bbb Z}_+$. Thus we get the following

\medskip

\begin{proposition} Let $\stackrel{\circ}{f}\in L_q^\infty$.
\begin{enumerate}
\item For all $\alpha>0$ one has a expansion in $D(U)_q^\prime$
$$B_{q,\alpha}\stackrel{\circ}{f}=\sum_{n \in{\Bbb Z}_+}q^{4 \alpha
n}B_q^{(n)}\stackrel{\circ}{f}.$$
\item For all $\psi \in D(U)_q$ one has the asymptotic expansion
$$\int \limits_{U_q}\left(B_{q,\alpha}\stackrel{\circ}{f}\right)\psi d
\nu {{\atop{\atop{\atop{\atop \displaystyle \sim}}}} \atop
{\scriptscriptstyle \alpha \to +\infty \atop}}\sum_{n=0}^\infty q^{4 \alpha
n}\int \limits_{U_q}\left(B_q^{(n)}\stackrel{\circ}{f}\right)\psi d
\nu.\eqno(5.2)$$
\end{enumerate}
Here $B_q^{(n)}:D(U)_q^\prime \to D(U)_q^\prime$ are polynomial functions of
the Laplace-Beltrami operator, given explicitly by (3.6).
\end{proposition}

\bigskip

\section{Covariant symbols}

 The notation $\widehat{z}$, $\widehat {z}^*$ in \cite{SSV1} stand for the
Toeplitz-Bergman operators with symbols $z$, $z^*$. Those are defined in the
graded vector space ${\Bbb C}[z]_{q,\alpha}$, with ${\rm
deg}(\widehat{z})=+1$, ${\rm deg}(\widehat{z}^*)=-1$. Hence for any matrix
$(a_{ij})_{i,j \in{\Bbb Z}_+}$ with numerical entries, series
$$\widehat{f}=\sum_{i,j \in{\Bbb Z}_+}a_{ij}\widehat{z}^i
\widehat{z}^{*j}\eqno(6.1)$$
converge in the topological vector space $\overline{F}_{q,\alpha}={\rm
Hom}_{\Bbb C}({\Bbb C}[z]_{q,\alpha},{\Bbb C}[[z]]_{q,\alpha})$.

\medskip

\begin{proposition} $\widehat{f}z^n=\displaystyle \sum \limits_{m \in{\Bbb
Z}_+}b_{mn}z^m$, $n \in{\Bbb Z}_+$,\\
with $b_{mn}=\displaystyle \sum
\limits_{j=0}^{\min(m,n)}\frac{\textstyle(q^{2n};q^{-2})_{n-j}}
{\textstyle(q^{4 \alpha+2n};q^{-2})_{n-j}}a_{m-j,n-j}$.
\end{proposition}

\smallskip

 {\bf Proof.} It suffices to apply the relations
$$\widehat{z}(z^m)=z^{m+1},\qquad \widehat{z}^*(z^m)=\left
\{\begin{array}{cll}\frac{\textstyle 1-q^{2m}}{\textstyle 1-q^{4
                    \alpha+2m}}\cdot z^{m-1} &,& m \ne 0 \\
                    0 &,&m=0 \end{array}\right.\eqno(6.2)$$
which were established in \cite[section 7]{SSV1} (see also \cite{KL}).\hfill
$\Box$

\medskip

\begin{corollary} For any linear operator $\widehat{f}\in
\overline{F}_{q,\alpha}$ there exists a unique decomposition (6.1).
\end{corollary}

\medskip

 {\sc Example 6.3.} Consider the linear operator $\widehat{f}_0:\,z^j
\mapsto \left \{\begin{array}{cll}1-q^{4 \alpha} &,& j=0 \\
                                  0 &,& j \ne 0 \end{array}\right.$, $j
\in{\Bbb Z}_+$. Prove that
$$\widehat{f}_0=(1-q^{4 \alpha})\sum_{k=0}^\infty \frac{(q^{-4
\alpha-2};q^2)_k}{(q^2;q^2)_k}q^{(4
\alpha+2)k}\widehat{z}^k \widehat{z}^{*k}.\eqno(6.3)$$

 Pass from the equality of operators to the equalities of their matricial
elements with respect to the base $\{z^n \}_{n \in{\Bbb Z}_+}$. Of course,
all the non-diagonal elements are zero. An identification of the diagonal
elements yields
$$\sum_{k=0}^j \frac{(q^{-4 \alpha-2};q^2)_k}{(q^2;q^2)_k}\cdot
\frac{(q^{2j};q^{-2})_k}{(q^{4 \alpha+2j};q^{-2})_k}\cdot q^{(4
\alpha+2)k}=\delta_{j0}.\eqno(6.4)$$
It suffices to consider the case $j>0$. Multiply (6.4) by
$\frac{\textstyle(q^{4 \alpha+2j};q^{-2})_j}{\textstyle(q^{2j};q^{-2})_j}$
to get
$$\sum_{k=0}^j \frac{(q^{-4 \alpha-2};q^2)_k}{(q^2;q^2)_k}\cdot
\frac{(q^{4 \alpha+2};q^2)_{j-k}}{(q^2;q^2)_{j-k}}\cdot
q^{(4 \alpha+2)k}=0.$$
That is,
$$\sum_{k+m=j} \frac{(q^{-4 \alpha-2};q^2)_k}{(q^2;q^2)_k}\cdot q^{(4
\alpha+2)k}\cdot \frac{(q^{4 \alpha+2};q^2)_m}{(q^2;q^2)_m}=0.$$
So, it remains to consider the q-binomial series (see \cite{GR}):
$$a(t)=\sum_{k \in{\Bbb Z}_+}\frac{(q^{-4
\alpha-2};q^2)_k}{(q^2;q^2)_k}\cdot q^{(4 \alpha+2)k}\cdot
t^k=\frac{(t;q^2)_\infty}{(q^{4 \alpha+2}t;q^2)_\infty},$$
$$b(t)=\sum_{m \in{\Bbb Z}_+}\frac{(q^{4 \alpha+2};q^2)_m}{(q^2;q^2)_m}\cdot
t^m=\frac{(q^{4 \alpha+2}t;q^2)_\infty}{(t;q^2)_\infty},$$
and to observe that $a(t)\cdot b(t)=1$.\hfill $\Box$

\medskip \stepcounter{theorem}

 It was noted in section 1 that $\widehat{f}_0$ is a Toeplitz-Bergman
operator with symbol $f_0$. This element generates the topological $U_q
\frak{sl}_2$-module $\overline{F}_{q,\alpha}$, as one can see from

\medskip

\begin{proposition} $U_q \frak{sl}_2 \widehat{f}_0=F_{q,\alpha}$.
\end{proposition}

\smallskip

 {\bf Proof.} Since for all $i,j,n \in{\Bbb Z}_+$,
$\widehat{z}^i \widehat{f}_0 \widehat{z}^{*n}:z^j \mapsto(1-q^{4
\alpha})\cdot \frac{\textstyle(q^2;q^2)_n}{\textstyle(q^{4 \alpha+2};q^2)_n}
\cdot \delta_{jn}z^i$, the linear operators $\{\widehat{z}^i \widehat{f}_0
\widehat{z}^{*n}\}_{i,n \in{\Bbb Z}_+}$ generate $F_{q,\alpha}$ as a vector
space. It remains to show that all those operators are in the $U_q
\frak{sl}_2$-module generated by $\widehat{f}_0$. For that, it suffices to
reproduce the proof of \cite[theorem 3.9]{SSV2}. One has only to alter the
notation for the generators (now they are $\widehat{z}$, $\widehat{z}^*$,
$\widehat{f}_0$), together with the constants in formulae which describe the
action of $X^\pm$ on $\widehat{f}_0$:
$$X^+\widehat{f}_0=c^\prime \widehat{z}\cdot \widehat{f}_0;\qquad
X^-\widehat{f}_0=c^{\prime \prime}\widehat{f}_0 \cdot \widehat{z}^*;\qquad
c^\prime,c^{\prime \prime}\ne 0.$$
These relations follow from proposition 1.2, \cite[proposition 3.8]{SSV2},
and
$${\Bbb C}\widehat{zf_0}={\Bbb C}\widehat{z}\widehat{f_0},\qquad{\Bbb
C}\widehat{f_0z^*}={\Bbb C}\widehat{f_0}\widehat{z^*}.$$
The latter relations can be deduced from
$${\rm Im}\,\widehat{zf_0}={\rm Im}\,\widehat{f_0z^*}={\Bbb C};\qquad{\rm
Ker}\,\widehat{f_0}={\rm Ker}\,\widehat{zf_0}={\Bbb C}^\perp.$$

\medskip

 It was shown in \cite[section 1]{SSV1} that for any $f \in D(U)_q^\prime$
there exists a unique decomposition $f=\sum \limits_{j,n \in{\Bbb
Z}_+}a_{jk}z^jz^{*k}$ similar to (6.1).

\medskip

 {\sc Example 6.5.} Prove that
$$(1-q^{4 \alpha})(1-zz^*)^{2 \alpha+1}=(1-q^{4 \alpha})\sum_{k \in{\Bbb
Z}_+}\frac{(q^{-(4 \alpha+2)};q^2)_k}{(q^2;q^2)_k}q^{(4
\alpha+2)k}z^kz^{*k}.\eqno(6.5)$$
Apply the operator $\widehat{T}$ to the both parts of (6.5) and identify the
matricial elements with respect to the base $\{z^m \}$ (it suffices
to consider the diagonal elements).

 Use the relations \footnote{These relations can be deduced from (6.2) via
passage to the limit as $\alpha \to \infty$,}
$$\widehat{T}(z)z^m=z^{m+1},\qquad \widehat{T}(z^*)z^m=\left
\{\begin{array}{cll}(1-q^{2m})z^{m-1} &,& m \ne 0 \\
                    0 &,& m=0 \end{array}\right.$$
to get
$$\sum_{k=0}^j \frac{(q^{-4 \alpha-2};q^2)_k}{(q^2;q^2)_k} \cdot
(q^{2j};q^{-2})_k \cdot q^{(4 \alpha+2)k}=q^{2j(2 \alpha+1)},$$
$$\sum_{k+m=j}\frac{(q^{-4 \alpha-2};q^2)_k}{(q^2;q^2)_k} \cdot q^{(4
\alpha+2)k}\cdot{1 \over(q^2;q^2)_m}=\frac{q^{2j(2
\alpha+1)}}{(q^2;q^2)_j}.$$
It remains to pass to the q-binomial decompositions (see \cite{GR}) in the
both sides of the obvious relation $a(t)b(t)=c(t)$, with
$$a(t)=\frac{(t;q^2)_\infty}{(q^{4 \alpha+2}t;q^2)_\infty};\qquad b(t)={1
\over (t;q^2)_\infty};\qquad c(t)={1 \over (q^{4
\alpha+2}t;q^2)_\infty}.\eqno \Box$$

\medskip \stepcounter{theorem}

\begin{proposition} The covariant symbol of the operator
$\widehat{f}=\displaystyle \sum \limits_{j,k \in{\Bbb
Z}_+}a_{jk}\widehat{z}^j \widehat{z}^{*k}$ is \\
$f=\displaystyle \sum \limits_{j,k \in{\Bbb Z}_+}a_{jk}z^jz^{*k}$.
\end{proposition}

\smallskip

 {\bf Proof.} Let $S_{q,\alpha}^\prime:\overline{F}_{q,\alpha}\to
D(U)_q^\prime$ be the map which takes a linear operator $\widehat{f}\in
\overline{F}_{q,\alpha}$ to its covariant symbol. We have to prove that this
map coincides with the map $S_{q,\alpha}^{\prime
\prime}:\overline{F}_{q,\alpha}\to D(U)_q^\prime$, given by
$S_{q,\alpha}^{\prime \prime}:\displaystyle \sum \limits_{j,k \in{\Bbb
Z}_+}a_{jk}\widehat{z}^j \widehat{z}^{*k} \mapsto \displaystyle \sum
\limits_{j,k \in{\Bbb Z}_+}a_{jk}z^jz^{*k}$. The linear operators
$S^\prime$, $S^{\prime \prime}$ are morphisms of $U_q \frak{sl}_2$-modules,
and the element $\widehat{f}_0$ generates the topological $U_q
\frak{sl}_2$-module $\overline{F}_{q,\alpha}$ by proposition 6.4. Thus it
suffices to obtain the relation $S^\prime \widehat{f}_0=S^{\prime
\prime}\widehat{f}_0$. It was shown in section 2 that
$S^\prime(\widehat{f}_0)=(1-q^{4 \alpha})(1-zz^*)^{2 \alpha+1}$. So it
remains to see that $S^{\prime \prime}(\widehat{f}_0)=(1-q^{4
\alpha})(1-zz^*)^{2 \alpha+1}$. This follows from (6.3), (6.5).\hfill $\Box$

\medskip

\begin{corollary} The map $\overline{F}_{q,\alpha}\to D(U)_q^\prime$ which
takes a linear operator to its covariant symbol is one-to-one.
\end{corollary}

\medskip

 To conclude, we give another illustration of corollary 6.2. Our immediate
purpose is to get the expansion $\widehat{z}^*\widehat{z}=\displaystyle
\sum \limits_{k \in{\Bbb Z}_+}c_k \widehat{z}^k \widehat{z}^{*k}$ and to find
a generating function $c(u)=\displaystyle \sum_{k \in{\Bbb Z}_+}c_ku^k$.

 By (6.2), the coefficients $c_k$ can be found from the system of equations
$$\sum_{k=0}^m c_k \frac{(q^{2m};q^{-2})_k}{(q^{4
\alpha+2m};q^{-2})_k}=\frac{1-q^{2(m+1)}}{1-q^{4 \alpha+2(m+1)}},\qquad m
\in{\Bbb Z}_+.\eqno(6.6)$$

 Apply an expansion of the right hand side of (6.6) as series:
$$\frac{1-q^{2(m+1)}}{1-q^{4 \alpha+2(m+1)}}=1+\sum_{j \in{\Bbb N}}(1-q^{-4
\alpha})q^{(2 \alpha+1+m)2j}.$$

 For a fixed $j \in{\Bbb N}$ consider the system of equations
$$\sum_{k=0}^m \gamma_k \frac{(q^{2m};q^{-2})_k}{(q^{4
\alpha+2m};q^{-2})_k}=q^{2mj},\qquad m \in{\Bbb Z}_+.\eqno(6.7)$$
Multiply (6.7) by $\frac{(q^{4 \alpha+2};q^2)_m}{(q^2;q^2)_m}$ and convert
it to the form
$$\sum_{i+k=m}\gamma_k \frac{(q^{4
\alpha+2};q^2)_i}{(q^2;q^2)_i}=q^{2mj}\frac{(q^{4
\alpha+2};q^2)_m}{(q^2;q^2)_m}.\eqno(6.8)$$
Introduce the generating functions
$$\alpha(u)=\sum_{m \in{\Bbb Z}_+}\frac{(q^{4
\alpha+2};q^2)_m}{(q^2;q^2)_m}q^{2mj}u^m=\frac{(q^{4
\alpha+2+2j}u;q^2)_\infty}{(q^{2j}u;q^2)_\infty},$$
$$\beta(u)=\sum_{i \in{\Bbb Z}_+}\frac{(q^{4
\alpha+2};q^2)_i}{(q^2;q^2)_i}u^i=\frac{(q^{4
\alpha+2}u;q^2)_\infty}{(u;q^2)_\infty}.$$

 It follows from(6.8) that
$$\gamma(u)\stackrel{\rm def}{=}\sum_{k \in{\Bbb
Z}_+}\gamma_ku^k={\alpha(u)\over \beta(u)}=\frac{(u;q^2)_j}{(q^{4
\alpha+2}u;q^2)_j}.$$
Turn back to the initial system (6.6) to obtain
$$c(u)=1+\sum_{j \in{\Bbb N}}(1-q^{-4 \alpha})q^{(2
\alpha+1)2j}\frac{(u;q^2)_j}{(q^{4 \alpha+2}u;q^2)_j}.\eqno(6.9)$$

\bigskip

\section{$*$ - Product}

 Let $A$ be an algebra over ${\Bbb C}$. Impose the notation
$${\Bbb C}[[q^{4 \alpha}]]=\left \{\sum_{n \in{\Bbb Z}_+}q^{4
\alpha}u_n|\:u_n \in{\Bbb C},\,n \in{\Bbb Z}_+\right \},$$
$$A[[q^{4 \alpha}]]=\left \{\sum_{n \in{\Bbb Z}_+}q^{4 \alpha}a_n|\:a_n \in
A \right \}$$
for the ring of formal series with complex coefficients and the ${\Bbb
C}[[q^{4 \alpha}]]$-algebra of formal series with coefficients from $A$.

 Our goal is to derive a new "distorted" multiplication in the ${\Bbb
C}[[q^{4 \alpha}]]$-algebra ${\rm Pol}({\Bbb C})_q[[q^{4 \alpha}]]$ from an
ordinary multiplication in the ${\Bbb C}[[q^{4 \alpha}]]$-algebra ${\rm
End}({\Bbb C}[z]_{q,\infty})[[q^{4 \alpha}]]$.

 The presence of the base $\{z^m \}_{m=0}^\infty$ in each vector space
${\Bbb C}[z]_{q,\alpha}$, ${\Bbb C}[z]_{q,\infty}$ allows one to "identify"
them via the isomorphisms $i_ \alpha:{\Bbb C}[z]_{q,\infty}\to{\Bbb
C}[z]_{q,\alpha}$;\ $i_ \alpha:z^m \mapsto z^m$, $m \in{\Bbb Z}_+$.

 Consider the linear operators $i_ \alpha^{-1} \widehat{z}^j
\widehat{z}^{*k}i_ \alpha$, $j,k \in{\Bbb Z}_+$ in ${\Bbb C}[z]_{q,\infty}$.
It follows from (6.2) that
$$i_ \alpha^{-1} \widehat{z}^j \widehat{z}^{*k}i_ \alpha:z^m \mapsto
\frac{(q^{2m};q^{-2})_k}{(q^{4 \alpha+2m};q^{-2})_k} \cdot z^{m-k+j}, \quad
m \in{\Bbb Z}_+.\eqno(7.1)$$

 From now on we shall identify the rational function ${\textstyle 1
\over \textstyle(q^{4 \alpha+2m};q^{-2})_k}$ of an indeterminate $t=q^{4
\alpha}$ with its q-binomial series (see \cite{GR})
$$\frac{(q^{4 \alpha+2m+2};q^2)_\infty}{(q^{4
\alpha+2m+2-2k};q^2)_\infty}=\sum_{n \in{\Bbb
Z}_+}\left(\frac{(q^{2k};q^2)_n}{(q^2;q^2)_n} \cdot
q^{2(m-k+1)n}\right)q^{4 \alpha n}.$$

 The construction of $*$-product will be done via the ${\Bbb C}[[q^{4
\alpha}]]$-linear map
$$Q:{\rm Pol}({\Bbb C})_q[[q^{4 \alpha}]]\to{\rm End}({\Bbb
C}[z]_{q,\infty})[[q^{4 \alpha}]]$$
defined as follows:
$$Q:\sum_{n \in{\Bbb Z}_+}q^{4 \alpha
n}\sum_{j,k=1}^{N(n)}a_{jk}^{(n)}z^jz^{*k} \mapsto \sum_{n \in{\Bbb
Z}_+}q^{4 \alpha n}\sum_{j,k=1}^{N(n)}a_{jk}^{(n)}i_
\alpha^{-1}(\widehat{z}^j \widehat{z}^{*k})i_ \alpha$$
for all numbers $a_{jk}^{(n)}\in{\Bbb C}$.

\medskip

\begin{lemma} The map $Q$ is injective.
\end{lemma}

\smallskip

 {\bf Proof.} In the case $Q$ has a non-trivial kernel, there should be for
some $j,k \in{\Bbb Z}_+$, $\displaystyle \sum \limits_{n \in{\Bbb Z}_+}c_n
\widehat{z}^{j+n}\widehat{z}^{*(k+n)}=0$, with $c_n \in{\Bbb C}[[q^{4
\alpha}]]$, $n \in{\Bbb Z}_+$, and $c_0 \ne 0$. An application of the
operator $\displaystyle \sum \limits_{n \in{\Bbb Z}_+}c_n
\widehat{z}^{j+n}\widehat{z}^{*(k+n)}$ to the vector $z^k$ yields $c_0 \cdot
\frac{(\textstyle q^{2k};q^{-2})_k}{\textstyle(q^{4 \alpha+2k};q^{-2})_k}
\cdot z^j=0$, which is a contradiction.\hfill $\Box$

\medskip

\begin{lemma} Let $j,k \in{\Bbb Z}_+$ and $\stackrel{\circ}{f}=z^{*j}z^k$.
The Toeplitz-Bergman operator $\widehat{f}$ with symbol
$\stackrel{\circ}{f}$ is $\widehat{z}^{*j}\widehat{z}^k$.
\end{lemma}

\smallskip

 {\bf Proof.} For all $\psi_1,\psi_2 \in H_{q,\alpha}^2$ one has
$$(\widehat{f}\psi_1,\psi_2)_{q,\alpha}=(P_{q,\alpha}(z^{*j}z^k \psi_1),
\psi_2)_{q,\alpha}=(z^{*j}z^k \psi_1, \psi_2)_{q,\alpha}=(z^k \psi_1,z^j
\psi_2)_{q,\alpha}=$$
$$=(\widehat{z}^k \psi_1,\widehat{z}^j
\psi_2)_{q,\alpha}=(\widehat{z}^{*j}\widehat{z}^k
\psi_1,\psi_2)_{q,\alpha}.\eqno \Box$$

\medskip

 The main result of this section is

\medskip

\begin{proposition} There exists a unique ${\Bbb C}[[q^{4
\alpha}]]$-bilinear map
$$*:{\rm Pol}({\Bbb C})_q[[q^{4 \alpha}]]\times{\rm
Pol}({\Bbb C})_q[[q^{4 \alpha}]]\to{\rm Pol}({\Bbb C})_q[[q^{4 \alpha}]]$$
such that $Q(f_1*f_2)=(Qf_1)\cdot(Qf_2)$ for all $f_1,f_2 \in{\rm Pol}({\Bbb
C})_q[[q^{4 \alpha}]]$.
\end{proposition}

\smallskip

 {\bf Proof.} The uniqueness follows from lemma 7.1. The existence of this
${\Bbb C}[[q^{4 \alpha}]]$-bilinear map will be established via verifying an
explicit formula (7.4). We start with considering the case $f_1=z^*$,
$f_2=z$.

 In section 6 a generating function $c(u)=\sum \limits_{k \in{\Bbb
Z}_+}c_ku^k$ for the coefficients of the expansion
$$\widehat{z}^*\widehat{z}=\sum_{k \in{\Bbb Z}_+}c_k \widehat{z}^k
\widehat{z}^{*k}\eqno(7.2)$$
was derived. Prove that
$$B_{q,\alpha}(z^*z)=\sum_{k \in{\Bbb Z}_+}c_k z^kz^{*k}.\eqno(7.3)$$
In fact, the distribution $B_{q,\alpha}(z^*z)$ coincides with the covariant
symbol of the Toeplitz-Bergman operator with symbol $z^*z$. This operator is
$\widehat{z}^*\widehat{z}$ by a virtue of corollary 7.2. Its covariant symbol
is $\sum \limits_{k \in{\Bbb Z}_+}c_k z^kz^{*k}$ due to proposition 6.6.

 It should be noted that $B_{q,\alpha}(z^*z)\in{\rm Pol}({\Bbb C})_q[[q^{4
\alpha}]]$. In fact, (6.9) implies
$$c(u)=c(u,q^{4 \alpha})=\sum_{n \in{\Bbb Z}_+}q^{4 \alpha n}\cdot
P_n(u),$$
with $P_n(u)$ being a polynomial of a degree at most $n+1$. Now our
statement in the case $f_1=z^*$, $f_2=z$ follows from (7.2) and (7.3):
$$z^**z=B_{q,\alpha}(z^*z).$$

 In a more general setting $f_1=z^{*m}$, $f_2=z^k$, $m,k \in{\Bbb Z}_+$, one
can use a similar argument. One has:
$$z^{*m}*z^k=B_{q,\alpha}(z^{*m}z^k).$$
The relations $Q(zf)=zQ(f)$, $Q(fz^*)=Q(f)z^*$, $f \in{\rm Pol}({\Bbb
C})_q[[q^{4 \alpha}]]$, allow one to consider even more general case of
$f_1,f_2 \in{\rm Pol}({\Bbb C})_q$:
$$z^iz^{*m}*z^kz^{*j}=z^iB_{q,\alpha}(z^{*m}z^k)z^{*j};\qquad i,j,k,m
\in{\Bbb Z}_+.\eqno(7.4)$$

 To complete the proof of proposition 7.3, it remains to define the
$*$-product of formal series:
$$\sum_{i \in{\Bbb Z}_+}q^{4 \alpha i}f_1^{(i)}*\sum_{j \in{\Bbb Z}_+}q^{4
\alpha j}f_2^{(j)}\stackrel{\rm def}{=}\sum_{n \in{\Bbb Z}_+}q^{4 \alpha
n}\left(\sum_{i+j=n}f_1^{(i)}*f_2^{(j)}\right),$$
with $f_1^{(i)},f_2^{(j)}\in{\rm Pol}({\Bbb C})_q$, $i,j \in{\Bbb
Z}_+$.\hfill $\Box$

\medskip

 {\sc Remark 7.4.} The polynomials $P_n(u)$, $n \in{\Bbb Z}_+$, could be
found without application of the {\sl explicit formula} for generating
function (6.9). In fact, if one sets up $c_k=\displaystyle \sum \limits_{n
\in{\Bbb Z}_+}q^{4 \alpha n}c_k^{(n)}$,
$$\widehat{z}^*\widehat{z}=\sum_{n \in{\Bbb Z}_+}q^{4 \alpha
n}\sum_{k=0}^{n+1}c_k^{(n)}\widehat{z}^k \widehat{z}^{*k}.\eqno(7.5)$$
The constants $c_k^{(n)}$ could be found from the relation (see \cite{SSV1})
$$\widehat{z}^*\widehat{z}=q^2 \widehat{z}\widehat{z}^*+1-q^2+q^{4
\alpha}\cdot \frac{1-q^2}{1-q^{4
\alpha}}\cdot(1-\widehat{z}\widehat{z}^*)(1-\widehat{z}^*\widehat{z}).$$
(For example, $P_0=c_0^{(1)}u+c_0^{(0)}=q^2u+1-q^2$.) This kind of
description for coefficients in (7.5) was used in \cite{SSV1}. We observe
that the $*$-product introduced here coincides with the $*$-product
considered in \cite{SSV1}.

\bigskip

\section{$*$ - Product and q-differential operators}

 The operators $\Box$, ${\textstyle \partial^{(l)}\over \textstyle \partial
z^*}$, ${\textstyle \partial^{(r)}\over \textstyle \partial
z^*}$, ${\textstyle \partial^{(l)}\over \textstyle \partial z}$,
${\textstyle \partial^{(r)}\over \textstyle \partial z}$ were introduced in
\cite{SSV1}.

\medskip

\begin{lemma} Let $\varphi,\psi$ be polynomials of one indeterminate. Then
$${\partial^{(r)}\over \partial
z^*}(\varphi(z^*)\psi(z))={\partial^{(r)}\varphi(z^*)\over \partial
z^*}\cdot \psi(q^2z).$$
\end{lemma}

\smallskip

 {\bf Proof.} Since $dz^* \cdot z=q^2z \cdot dz^*$, one has
$${\partial^{(r)}\over \partial
z^*}(\varphi(z^*)\psi(z))\cdot dz^*=\overline{\partial}(\varphi(z^*)\psi(z))=
(\overline{\partial}\varphi(z^*))\psi(z)=$$
$$={\partial^{(r)}\varphi(z^*)\over \partial z^*}\cdot dz^*\cdot
\psi(z)={\partial^{(r)}\varphi(z^*)\over \partial z^*}\cdot \psi(q^2z)\cdot
dz^*.\eqno \Box$$

\medskip

\begin{lemma} For all $\psi(z)\in{\Bbb C}[z]_q$,
${\textstyle \partial^{(r)}\psi(z)\over \textstyle \partial z}={\textstyle
\partial^{(l)}\psi(q^2z)\over \textstyle \partial z}$.
\end{lemma}

\smallskip

 {\bf Proof.} Since $dz \cdot z=q^2z \cdot dz$, one has
$${\partial^{(r)}\psi(z)\over \partial z}\cdot dz=\partial \psi=dz
\cdot{\partial^{(l)}\psi(z)\over \partial z}={\partial^{(l)}\psi(q^2z)\over
\partial z}\cdot dz.\eqno \Box$$

\medskip

\begin{proposition} Let $f_1,f_2$ be polynomials of one indeterminate. Then
$$\Box(f_2(z^*)f_1(z))=q^{2}\cdot{\partial^{(r)}f_2 \over \partial
z^*}\cdot(1-zz^*)^2 \cdot{\partial^{(l)}f_1 \over \partial z}.\eqno(8.1)$$
\end{proposition}

\smallskip

 {\bf Proof.} It follows from \cite[corollary 2.9]{SSV4} that
$$\Box(f_2(z^*)f_1(z))=q^2 \left({\partial^{(r)}\over \partial
z^*}{\partial^{(r)} \over \partial z}(f_2(z^*)f_1(z)))\right)(1-zz^*)^2=$$
$$=q^{-2}\left({\partial^{(r)}\over \partial
z^*}\left(f_2(z^*){\partial^{(r)}f_1(z)\over \partial
z}\right)\right)(1-z^*z)^2.$$
Apply lemmas 8.1, 8.2 to conclude that
$$\Box(f_2(z^*)f_1(z))=q^{-2}{\partial^{(r)}f_2(z^*)\over \partial
z^*}\cdot{\partial^{(l)}f_1(q^4z)\over \partial z}(1-z^*z)^2.$$
It remains to apply the commutation relation
$z(1-z^*z)^2=q^{-4}(1-z^*z)^2z$.\hfill $\Box$

\medskip

 Remind the notation from \cite{SSV1}:
$$\stackrel{\sim}{\Box}=q^{-2}(1-(1+q^{-2})z^*\otimes
z+q^{-2}z^{*2}\otimes z^2)\cdot{\partial^{(r)}\over \partial
z^*}\otimes{\partial^{(l)}\over \partial z},$$
$m:{\rm Pol}({\Bbb C})_q \otimes{\rm Pol}({\Bbb C})_q \to{\rm Pol}({\Bbb
C})_q$, $m:\psi_1 \otimes \psi_2 \mapsto :\psi_1 \psi_2$.

 Now we are in a position to prove \cite[theorem 7.3]{SSV1}.

\medskip

\begin{theorem} For all $f_1,f_2 \in{\rm Pol}({\Bbb C})_q$
$$f_1*f_2=(1-q^{4 \alpha}) \cdot \sum_{j \in{\Bbb Z}_+}q^{4
\alpha \cdot j}m(p_j(\stackrel{\sim}{\Box})f_1 \otimes
f_2),$$
with $p_j$, $j \in{\Bbb Z}_+$, being the polynomials determined by (3.3).
\end{theorem}

\smallskip

 {\bf Proof.} With $f_1, f_2, f_3, f_4 \in{\Bbb C}[z]_q$, one can deduce
from the results of section 7 that
$$(f_1(z)f_2(z)^*)*(f_3(z)f_4(z)^*)=
f_1(z)B_{q,\alpha}(f_2(z)^*f_3(z))f_4(z)^*.\eqno(8.1)$$
An application of the results of section 5 to the bounded function
$f_2(z)^*f_3(z)$ yields:
$$B_{q,\alpha}(f_2(z)^*f_3(z)){{\atop \textstyle \sim}\atop{\alpha \to
\infty \atop}}(1-q^{4 \alpha})\sum_{j \in{\Bbb Z}_+}q^{4
\alpha \cdot j}p_j(\Box)(f_2(z)^*f_3(z)).$$
It remains to apply proposition 8.3 and the definition of
$\stackrel{\sim}{\Box}$.\hfill $\Box$

\medskip

 {\sc Remark 8.5.} One can observe from proposition 2.5 that (8.1) is a
q-analogue of relation (4.7) from \cite{B}.

\bigskip

\section*{Appendix. Overflowing vector systems}

 Unlike the main text where $\alpha$ was allowed to be an arbitrary positive
number, let us assume now $\alpha \in{1 \over 2}{\Bbb N}$.

 Remind the notation $\widetilde{X}$ for the quantum principal homogeneous
space, and $i:D(U)_q^\prime \hookrightarrow D(\widetilde{X})_q$ for the
canonical embedding of distribution spaces (see \cite{SSV3}).

 Consider the embedding of vector spaces
$$i_ \alpha:{\rm Pol}({\Bbb C})_q \hookrightarrow D(\widetilde{X})_q^\prime;
\qquad i_ \alpha:f \mapsto i(f)\cdot t_{12}^{-2 \alpha-1}.$$
Equip ${\rm Pol}({\Bbb C})_q$ with a new $U_q \frak{sl}_2$-module structure
given by $i_ \alpha \xi f=\xi i_ \alpha f$ for all $f \in{\rm Pol}({\Bbb
C})_q$, $\xi \in U_q \frak{sl}_2$. Denote this $U_q \frak{sl}_2$-module by
${\rm Pol}({\Bbb C})_{q,\alpha}$. There exists an embedding ${\Bbb
C}[z]_{q,\alpha}\hookrightarrow{\rm Pol}({\Bbb C})_{q,\alpha}$.

 The results of \cite[section 6]{SSV3} imply

\medskip

 {\bf Proposition A.1.} {\it The linear map $D(\widetilde{X})_q \to{\Bbb
C}[z]_{q,\alpha}$ given by
$$\psi \mapsto \int \limits_{\widetilde{X}_q}\tau_{12}^{*(-2
\alpha-1)}\cdot(z \zeta^*;q^2)_{2 \alpha+1}^{-1} \cdot \psi d \nu,$$
is a morphism of $U_q \frak{sl}_2$-modules.}

\medskip

 Proposition A.1 allows one to treat the function $\tau_{12}^{*(-2
\alpha-1)}\cdot(z \zeta^*;q^2)_{2 \alpha+1}^{-1}$ as a q-analogue of a
coherent state in the sense of Perelomov \cite{P}.

\medskip

 {\bf Corollary A.2.} {\it For all $\psi \in D(U)_q$
$$P_{q,\alpha}\psi(z)=\int \limits_{U_q}(z \zeta^*;q^2)_{2
\alpha+1}^{-1}\psi(\zeta)d \nu_ \alpha(\zeta).\eqno(A.1)$$}

\smallskip

 {\bf Proof.} Consider the integral operator
$$P:D(U)_q \to{\Bbb C}[z]_{q,\alpha};\qquad P:\psi(z)\mapsto \int
\limits_{U_q}(z \zeta^*;q^2)_{2 \alpha+1}^{-1}\psi(\zeta)d \nu_ \alpha.$$
It is a morphism of $U_q \frak{sl}_2$-modules, as one can deduce from
proposition A.1. The orthoprojection $P_{q,\alpha}$ is also a morphism of
$U_q \frak{sl}_2$-modules, due to the invariance of the scalar product in
$H_{q,\alpha}^2$. It remains to use the relations $Pf_0=1-q^{4 \alpha}$,
$P_{q,\alpha}f_0=1-q^{4 \alpha}$, together with the fact that $f_0$
generates the $U_q \frak{sl}_2$-module $D(U)_q$ (see \cite{SSV2}).\hfill
$\Box$

\medskip

 {\sc Remark A.3.} (A.1) means that the distribution $(z \zeta^*;q^2)_{2
\alpha+1}^{-1}$ is a reproducing kernel.

\medskip

 Let us find the kernel of the integral operator
$P_{q,\alpha}\stackrel{\circ}{f}P_{q,\alpha}$. For $\stackrel{\circ}{f},\psi
\in D(U)_q$ one has by corollary A.2
$$P_{q,\alpha}\stackrel{\circ}{f}P_{q,\alpha}:\psi(z)\mapsto \int
\limits_{U_q}K_q(\stackrel{\circ}{f};z,z^\prime)\psi(z^\prime)d \nu_
\alpha(z^\prime),$$
with
$$K_q(\stackrel{\circ}{f};z,z^\prime)=
\int \limits_{U_q}(z \zeta^*;q^2)_{2
\alpha+1}^{-1}\stackrel{\circ}{f}(\zeta)\cdot(\zeta z^{\prime*};q^2)_{2
\alpha+1}^{-1}d \nu_ \alpha(\zeta).$$
Now an application of lemma 1.4 yields
$$K_q(\stackrel{\circ}{f};z,z^\prime)={1-q^{4 \alpha}\over 1-q^2}
\int \limits_{U_q}(z \zeta^*;q^2)_{2
\alpha+1}^{-1}\stackrel{\circ}{f}(\zeta)(\zeta z^{\prime*};q^2)_{2
\alpha+1}^{-1}(1-\zeta \zeta^*)^{2 \alpha+1}d \nu(\zeta)=$$
$$={1-q^{4 \alpha}\over 1-q^2}
\int \limits_{U_q}(\zeta z^{\prime*};q^2)_{2 \alpha+1}^{-1}(1-\zeta
\zeta^*)^{2 \alpha}(z \zeta^*;q^2)_{2 \alpha+1}^{-1}\cdot
\stackrel{\circ}{f}(1-\zeta \zeta^*)d \nu(\zeta)=$$
$$={1-q^{4 \alpha}\over 1-q^2}
\int \limits_{U_q}(1-\zeta \zeta^*)(\zeta z^{\prime*};q^2)_{2
\alpha+1}^{-1}(1-\zeta \zeta^*)^{2 \alpha}(z \zeta^*;q^2)_{2
\alpha+1}^{-1}\stackrel{\circ}{f}d \nu(\zeta).$$

 Finally, use the relation $(1-\zeta \zeta^*)\zeta=q^2 \zeta(1-\zeta
\zeta^*)$ to obtain

\medskip

 {\bf Proposition A.4.} {\it $P_{q,\alpha}\stackrel{\circ}{f}P_{q,\alpha}$
is an integral operator:
$$P_{q,\alpha}\stackrel{\circ}{f}P_{q,\alpha}\psi(z)=\int
\limits_{U_q}K_q(\stackrel{\circ}{f};z,z^\prime)\psi(z^\prime)d \nu_
\alpha(z^\prime),$$
whose kernel is given by
$$K_q(\stackrel{\circ}{f};z,z^\prime)={1-q^{4 \alpha}\over 1-q^2}\int
\limits_{U_q}k_ \zeta(q^2z^\prime)^* \cdot k_
\zeta(z)\stackrel{\circ}{f}(\zeta)d \nu(\zeta),$$
with $k_ \zeta(z)=(1-\zeta \zeta^*)^{\alpha+1/2}\cdot(\zeta^*z;q^2)_{2
\alpha+1}^{-1}\in D(U \times U)_q^\prime$.}

\medskip

 Proposition A.4 allows one to treat the distribution $k_ \zeta(z)$ as a
q-analogue of an overflowing vector system.

\end{document}